\documentclass[final,3p]{elsarticle}

\usepackage{amssymb}
\usepackage{amsthm}
\usepackage{amsmath}
\usepackage[table]{xcolor}

\renewcommand{\vec}[1]{{\mathchoice
                     {\mbox{\boldmath$\displaystyle{#1}$}}
                     {\mbox{\boldmath$\textstyle{#1}$}}
                     {\mbox{\boldmath$\scriptstyle{#1}$}}
                     {\mbox{\boldmath$\scriptscriptstyle{#1}$}}}}

\newcommand{\mat}[1]{\vec{\mathrm{#1}}}

\newcommand{\ta}{\boldsymbol{\theta}}

\newcommand{\tab}{\boldsymbol{\overline{\theta}}}
\newcommand{\tat}{\boldsymbol{\widetilde{\theta}}}
\newcommand{\tah}{\boldsymbol{\widehat{\theta}}}

\newcommand{\ff}{\boldsymbol{f}}

\newcommand{\m}{\vec{m}}
\newcommand{\ms}{\vec{m}^*}

\newcommand{\D}{\mat{D}}
\renewcommand{\H}{\mat{H}}
\newcommand{\B}{\mat{B}}

\newcommand{\R}{\mathbb{R}}

\newcommand{\pdiff}[2]{\frac{\partial {#1}}{\partial{#2}}}
\newcommand{\diff}[2]{\frac{d{#1}}{d{#2}}}

\newcommand{\eps}{\varepsilon}
\newcommand{\Newton}{\text{Newton}}

\journal{}

\begin{document}

\begin{frontmatter}



\title{A new perspective on parameter study of optimization problems}


\author[a]{Alen Alexanderian}
\author[b]{Joseph Hart \footnote{Corresponding author email address: joshart@sandia.gov.}}
\author[a]{Mason Stevens}
\address[a]{Department of Mathematics, North Carolina State University, Raleigh, NC}
\address[b]{Department of Scientific Machine Learning, Sandia National Laboratories, Albuquerque, NM}

\begin{abstract}
We provide a new perspective on the study of 
parameterized optimization problems.  Our approach combines methods for 
post-optimal sensitivity analysis and ordinary differential equations to quantify the uncertainty 
in the minimizer due to uncertain parameters in the optimization problem. 
We illustrate the
proposed approach with a simple analytic example and an inverse problem
governed by an advection diffusion equation.
\end{abstract}




\end{frontmatter}

\section{Introduction}


A common class of problems in the sciences and engineering involves solving
optimization problems constrained by differential equations. Examples
include inverse problems and optimal design or control problems.
To illustrate, we consider the following advection diffusion equation:
\begin{equation}\label{equ:CD}
\begin{aligned}
-\kappa u'' + v u' &= s \quad & \text{in } (0,1), \\
\kappa u' &= \alpha u     \quad & \text{on }     x = 0,\\ 
\kappa u' &= -\alpha u  \quad & \text{on } x=1.\\
\end{aligned}
\end{equation}
Here $u(x)$ is the temperature at a point $x \in [0,1]$, $\kappa$ is the diffusion coefficient, 
$v$ is wind velocity, $\alpha$ models a heat transfer coefficient, and $s$ is a source term defined by 
\begin{equation}\label{equ:f}
s(x) = a\exp \left(-200(x-c)^2 \right),
\end{equation}
which models a localized source with $a$ and $c$ indicating its magnitude and location.

Suppose we have measurements of the temperature throughout the domain and we seek to use
this information to estimate the parameter vector $\vec m = [\begin{matrix} \kappa & v \end{matrix}]^\top \in \R^2$.
This involves solving an optimization problem of the form,
\begin{equation}\label{equ:optim_CD}
\min_{\vec m} J(\vec m) := \frac{1}{2} \int_0^1 \big(u(x)-u^\text{obs}(x)\big)^2 \, dx + \frac{\beta}{2} \| \vec m - \vec{m}^0 \|_2^2,
\end{equation}
where $u$ is the solution of~\eqref{equ:CD} (which depends on $\vec m$), $u^\text{obs}$ is 
the observed temperature across the domain, $\vec{m}^0$ is a prior estimate 
of $\vec m$, and  $\beta > 0$ is a regularization parameter. 
Note that the objective function $J$ also depends on the vector of model
parameters $\vec \theta = [\begin{matrix}  a & c & \alpha \end{matrix}]^\top
\in \R^3$ that parameterize the volume source term and the heat transfer.  In
practice, these parameters might not be known exactly.  Thus, it is imperative
to understand how the uncertainty in these parameters affects the estimated
parameter $\ms$ obtained by solving~\eqref{equ:optim_CD}. In the present work,
we propose an approach for analyzing such \emph{parameterized} optimization
problems, without the need for repeated solutions of the optimization problem.  

The study of parameterized optimization problems can be found in the early 
works~\cite{fiacco_mccormick_1968,fiacco_1976} followed by more advanced
developments in~\cite{bonnans_1992,shapiro_SIAM_review,bonnans_shapiro_book}.
These works, and references therein, provide extensive theory concerning the
differentiability of optimal solutions with respect to parameter perturbations.
This field has assumed various names as it arises in different parts of the
literature. Herein we refer to it as post-optimality sensitivity analysis as it
provides a local sensitivity study of the optimal solution. Developments
from~\cite{Griesse_part_1,Griesse_part_2,Griesse_Thesis,Brandes06} extended the
use of post-optimality sensitivities to optimization problems constrained by
partial differential equations. Recent work
from~\cite{HartvanBloemenWaanders20,SunseriHartEtAl20,rgsvd_saibaba} has
focused on making this sensitivity analysis scalable for high-dimensional
parameter spaces and extending its use for various classes of parametric
uncertainty. However, post-optimality sensitivity analysis is local in the
sense that it is only valid in a neighborhood of a nominal parameter value.
This article borrows concepts from the post-optimal sensitivity analysis
literature and couples them with a time stepping algorithm to move through the
parameter space to perform efficient global analysis. 

We detail the mathematical setup of the parameterized optimization problems
under study in Section~\ref{sec:prelim}. Our proposed approach is presented
in Section~\ref{sec:method}.  We illustrate the effectiveness of our approach
in Section~\ref{sec:numerics}, for an analytic test problem as well as the
inverse advection diffusion problem discussed above.  Concluding remarks
are given in Section~\ref{sec:conc}.

\section{Preliminaries}\label{sec:prelim}
In this section, we lay out the mathematical setup of 
the optimization problems under study.  
Let $U$ be a compact subset of $\R^d$ and consider parameterized optimization problems of the form
\begin{equation}\label{equ:optim}
    \min_{\m \in U} J(\m, \ta), 
\end{equation}
where $\ta$ is a vector of parameters. These parameters are fixed when
solving the optimization problem, but in practice are uncertain and can be modeled
as random variables. We
assume that $\ta$ belongs to a compact set $\Theta \subset \R^p$ and $\bar\ta
\in \Theta$ is a nominal parameter vector. Let $\bar\m^*$ be a minimizer of
$J(\m, \bar\ta)$ and $U_0 \subset U$ be an open set that contains $\bar\m^*$.
To facilitate our parametric study of~\eqref{equ:optim}, we assume that
\begin{enumerate}
\item for each $\ta \in \Theta$, there exists a unique minimizer $\ms(\ta)$ in 
$U_0$; and 
\item $J$ is twice continuously differentiable in $\vec m$ and $\ta$.
\end{enumerate}
%
It follows that $\ms(\ta)$ is a differentiable function on $\Theta$.

Verifying the above assumptions for an optimization problem arising from
science or engineering applications would be difficult in general. However,
these assumptions have an intuitive interpretation. Namely, we consider
well-behaved optimization problems whose minimizers are unique in a suitable
region around nominal parameters $\bar\ta$ and minimizer $\bar\m^*$. 
Such assumptions are reasonable in many optimization problems 
arising in physical applications. Also, in practice, the set $\Theta$ will be a user-specified
region around the nominal parameter vector $\tab$. Specifically, we consider  
the typical situation where  
$\Theta$ is defined as  
\begin{equation}\label{equ:Theta}
    \Theta = 
    [\bar\theta_1 - \eps_1, \bar\theta_1 + \eps_1] \times 
    [\bar\theta_2 - \eps_2, \bar\theta_2 + \eps_2] \times 
    \cdots \times
    [\bar\theta_p - \eps_p, \bar\theta_p + \eps_p],
\end{equation}
where $\eps_k$ is some percentage of the corresponding nominal value
$\bar\theta_k$, $k = 1, \ldots, p$.  The $\eps_k$'s indicate the level of uncertainty in
the physical parameters in the model.  It is common that the parameters $\theta_k$ are assumed to 
be uniformly distributed random variables on the respective intervals, but more general distributions on compact sets are admissible.

\begin{figure}\centering
\includegraphics[width=0.35\textwidth]{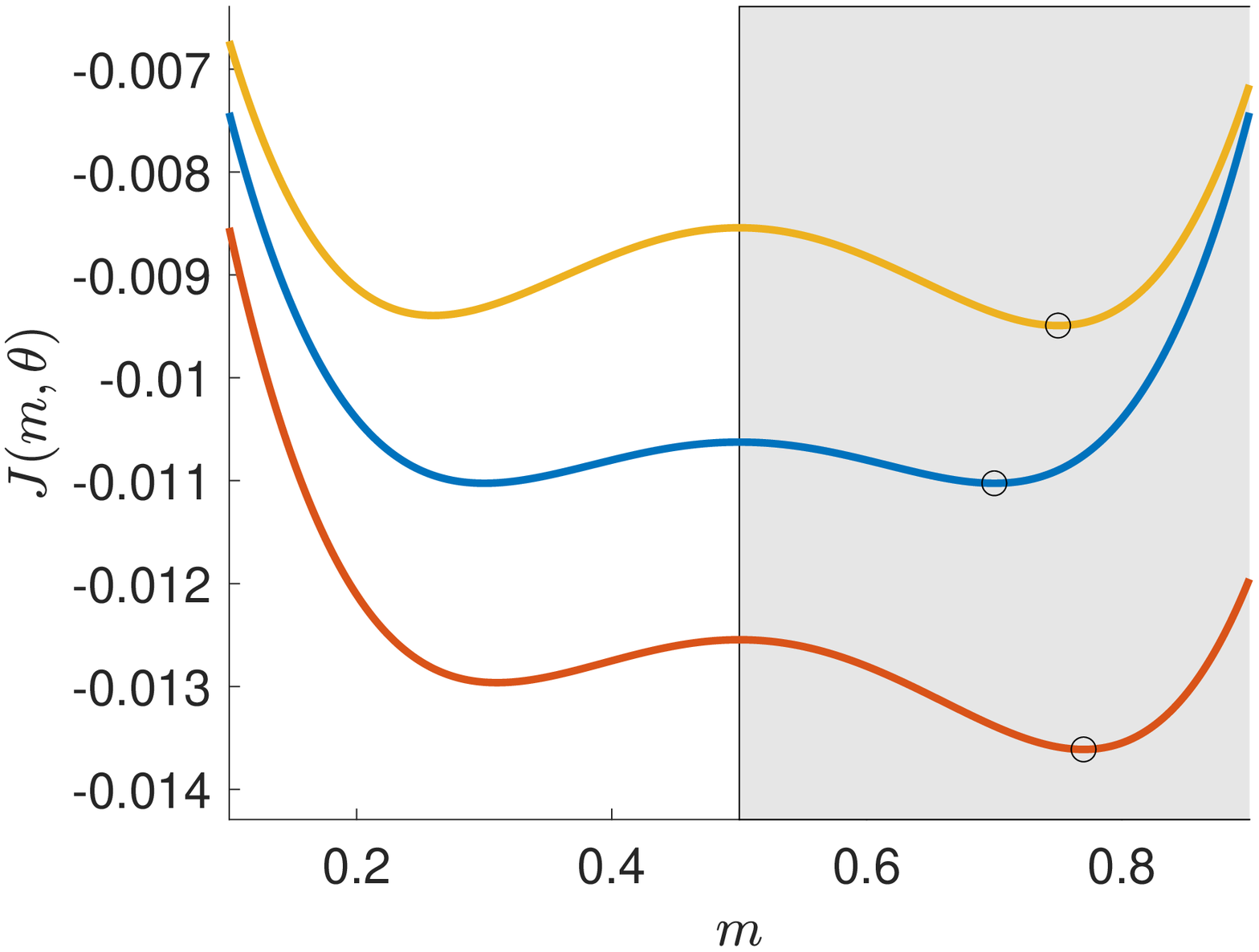}
\includegraphics[width=0.35\textwidth]{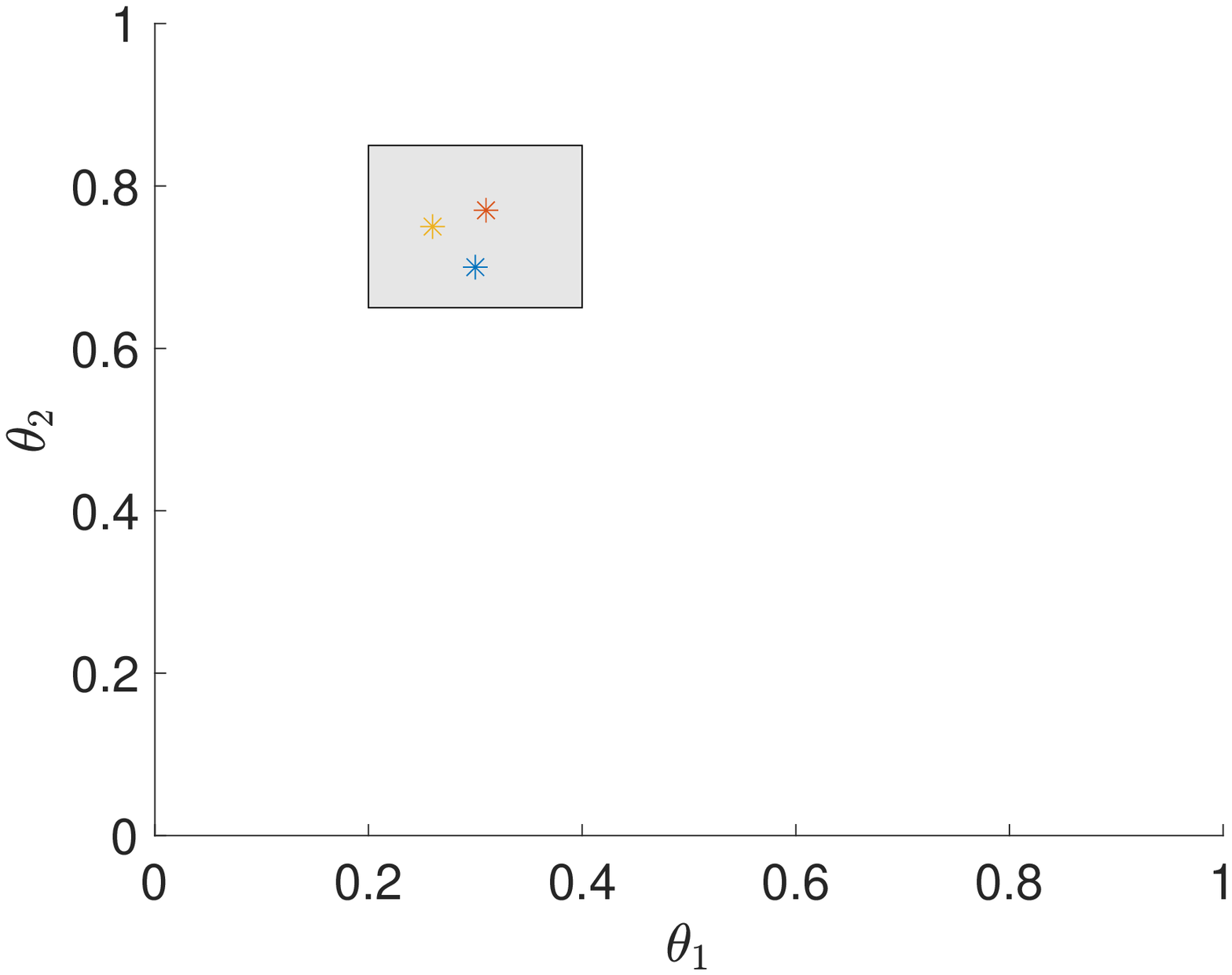}
\caption{Illustration of the sets $U_0=[0.5,1.0]$ shaded in the left panel and $\Theta=[0.2,0.4] \times [0.65,0.85]$ shaded in the right panel. Each curve in the left panel is $J(m,\ta)$, defined by~\eqref{equ:J_illustration}, evaluated for a $\ta$ sample depicted in the right panel by $*$. The unique local minimizers in $U_0$ are denoted by open circles in the left panel.}
\label{fig:local_min}
\end{figure}

To illustrate these concepts, consider the objective function 
\begin{eqnarray}
\label{equ:J_illustration}
J(m,\ta) = \int (m-\theta_1)(m-0.5)(m-\theta_2) dm,
\end{eqnarray}
which, for $\ta \in \{\ta \in \R^2 : 0 < \theta_1 <0.5 \text{ and } 0.5 < \theta_2 <1\}$, will
have two local minima at $\theta_1$ and $\theta_2$. Letting 
$U=[0,1]$ and $\Theta=[0.2,0.4] \times [0.65,0.85]$, we may
take $U_0 = (0.5,1.0)$ to restrict our analysis to the local minima at
$\theta_2$. This is depicted in Figure~\ref{fig:local_min} to demonstrate how
the choice of $U_0$ ensures minimizer uniqueness needed to enable our parameter
study.

\section{Method}\label{sec:method}
In this section, we outline an approach for approximating 
$\ms(\tat)$, where $\tat$ is a generic element of $\Theta$. 

\subsection{An initial value problem for $\ms(\tat)$}
To study how the optimal solution changes with $\ta$ we 
compute the Jacobian of $\ms$ with respect to $\ta$. 
This Jacobian, denoted by $\D$, can be 
computed by differentiating through first order optimality condition
\[
    \frac{\partial J}{\partial m}(\ms(\ta), \ta) = 0,
\]
implicitly. It follows from the Implicit Function Theorem that
\begin{equation*}
		\D(\ms(\ta),\ta) = -\H(\ms(\ta),\ta)^{-1}
                                 \B(\ms(\ta),\ta),
\end{equation*}
where
\begin{equation*}
		\H(\ms(\ta),\ta) = \frac{\partial^2 J}{\partial \m^2}\Big|_{(\ms(\ta),\ta)} \qquad \text{and} \qquad
		\B(\ms(\ta),\ta) = \frac{\partial^2 J}{\partial \m \partial \vec\theta}
                                 \Big|_{(\ms(\ta),\ta)}. 
\end{equation*}
The Jacobian $\D$ is known as the post-optimality sensitivity operator.

For a given $\tat \in \Theta$, we begin by considering the points on the
line-segment joining $\tab$ to $\tat$,
\begin{align}\label{equ:def_theta_t}
\vec\theta(t) = \tab + t\,( \tat - \tab), \quad t \in [0, 1].
\end{align}
Based on our assumptions on $J$ and $\ms$, we can consider $\ms$ as a differentiable 
function of $t$, $\ms(t) \equiv \ms(\vec\theta(t))$.
Taking the derivative of $\ms$ with respect to
$t$ and applying the Chain Rule gives
	\begin{align}\label{equ:dmdt}
		\diff{\ms}{t} = \pdiff{\ms}{\vec\theta} \cdot \diff{\vec\theta}{t} = \D(\ms(\ta(t)),\vec\theta(t))  ( \tat - \tab).
	\end{align}
For notational convenience, we define 
\begin{equation}\label{equ:fdef}
   \ff(t, \ms) = -\H(\ms,\ta(t))^{-1} \B(\ms,\ta(t))(\tat - \tab).
\end{equation}
Then we may determine $\ms(t)$, for $t \in [0, 1]$, by solving the following initial value problem (IVP) 
\begin{equation}\label{equ:IVP}
   \begin{aligned}
   \diff{\ms}{t} &= \ff(t, \ms), \\
\ms(0) &= \ms(\tab).
\end{aligned}
\end{equation}
To specify the initial condition, the optimization problem~\eqref{equ:optim}
needs to be solved.  This is the only solution of the optimization problem
required in our approach. For each parameter $\tat$, we will solve the IVP~\eqref{equ:IVP} up to $t=1$ to determine the corresponding local minimum $\ms(1) \equiv \ms(\ta(1)) = \ms(\tat)$. 
The right hand side function $\ff$~\eqref{equ:fdef} is the post optimality sensitivity operator acting on $\tat-\tab$ and hence the IVP depends on the parameter sample $\tat$. 

\subsection{Time-stepping to approximate $\ms(\tat)$}
We can apply common numerical methods to solve the IVP~\eqref{equ:IVP}. 
In this work, we study the use of forward Euler. Specifically, 
let $h = 1/N$ be a step-size and 
$t_n = n h$, $n = 0, \ldots, N$. The forward Euler discretization 
of~\eqref{equ:IVP} is
\begin{align}\label{equ:euler}
   \ms_{n+1} = \ms_n + h \, \ff(t_n, \ms_n), 
   \quad
  n=0,1,\ldots,N-1,
	\end{align}
where $\ms_0 = \ms(0)$ and $\ms_n \approx \ms(t_n)$, $n = 1, \ldots, N$. Finally, the approximation to 
$\ms(\tat)$ is given by 
\begin{equation}\label{equ:approx}
    \ms(\tat) \equiv \ms(t_N) \approx \ms_N.
\end{equation}

 Notice that~\eqref{equ:euler} resembles Newton's method for optimization. In particular, computing $\ms(\tat)$ via Newton's method iterates with search directions of the form
 \begin{eqnarray*}
 -\H(\vec{m}_n^{\Newton},\tat)^{-1} \vec{g}(\vec{m}_n^{\Newton},\tat),
 \end{eqnarray*}
 where $\vec{m}_n^{\Newton}$ denotes the $n^{th}$ Newton iterate and $\vec{g}(\vec{m}_n^{\Newton},\tat)$ denotes the gradient of $J$ with respect to $\vec{m}$, evaluated at $\vec{m}_n^{\Newton}$. On the other hand, time marching via~\eqref{equ:euler} has search directions of the form
  \begin{eqnarray*}
 -\H(\ms_n,\ta(t_n))^{-1} \B(\ms_n,\ta(t_n))h(\tat-\tab) \approx 
-\H(\ms_n,\ta(t_n))^{-1} \big(\vec{g}(\ms_n,\ta(t_{n+1})) - \vec{g}(\ms_n,\ta(t_{n}))\big),
 \end{eqnarray*}
where the latter approximation follows since $\B$ is the derivative of
$\vec{g}$ with respect to $\ta$. Note that 
$\vec{g}(\ms_n,\ta(t_{n}))\big)$ is expected to be small since it is the gradient evaluated at an approximate minimizer. Time marching via~\eqref{equ:euler} takes a
size $h$ perturbation of $\ta$ and updates the solution via a Newton like step
whereas resolving the optimization problem takes the full parameter step
$\tat-\tah$ and then employs Newton iterations to update $\vec{m}$.
Note also that in resolving the optimization problem via Newton's method, employing a line
search algorithm such as Armijo's method~\cite{NocedalWright99} is typically necessary. 
An alternate point of view regarding the time stepping~\eqref{equ:euler} is 
performing continuation on the parameters $\ta_0,\ta_1,\dots,\ta_N$ and using an
approximate Newton step to update $\vec{m}$ after each parameter perturbation.

A benefit of solving~\eqref{equ:IVP} is that computing $\ff(t_n,
\ms_n)=\D(\ms_n,\ta_n)(\tat-\tab)$ at each time step provides post-optimality sensitivity information in the direction $\tat-\tab$. This gives additional insight about how the minimizer depends on the uncertain parameters.

\subsection{Computational considerations}
Computing~\eqref{equ:fdef} requires access to the Hessian $\H$ and the matrix $\B$ of mixed
second order partial derivatives.  For optimization problems governed by partial
differential equations, such derivative information can be obtained efficiently
using adjoint state methods; see e.g.,~\cite{Gunzburger03}.  Specifically, in that context,
one obtains adjoint based expressions for applying $\H$ and $\B$ to 
vectors~\cite{SunseriHartEtAl20}. In large-scale computations, the
inverse Hessian apply is computed by performing a linear solve using
the Conjugate-Gradient method, which only requires 
applications of the Hessian on vectors. 
An alternative approach for obtaining the required derivatives is automatic
differentiation. Simple finite-difference approaches might be applicable as
well if the gradient can be computed exactly and differenced to approximate the
second derivatives. 

It is also possible to use methods other than forward Euler to numerically
solve the IVP~\eqref{equ:IVP}. Generally, implicit methods would be very
challenging to implement for the present IVP, due the requirement of a
nonlinear solve in each step. On the other hand, higher order explicit
Runge--Kutta methods or predictor corrector methods will be straightforward to
implement.  However, we caution that the faster convergence might come at a
cost of making the time-stepping more expensive than resolving the optimization
problem for different realizations of $\tat$. A thorough investigation of the
time-stepping approaches that are tractable for~\eqref{equ:IVP} and analysis of
the related computational cost is beyond scope of the present study and will be
pursued in future work.

\section{Numerical examples}\label{sec:numerics}
\subsection{A one-dimensional example}
\begin{figure}[ht]\centering
\includegraphics[width=0.32\textwidth]{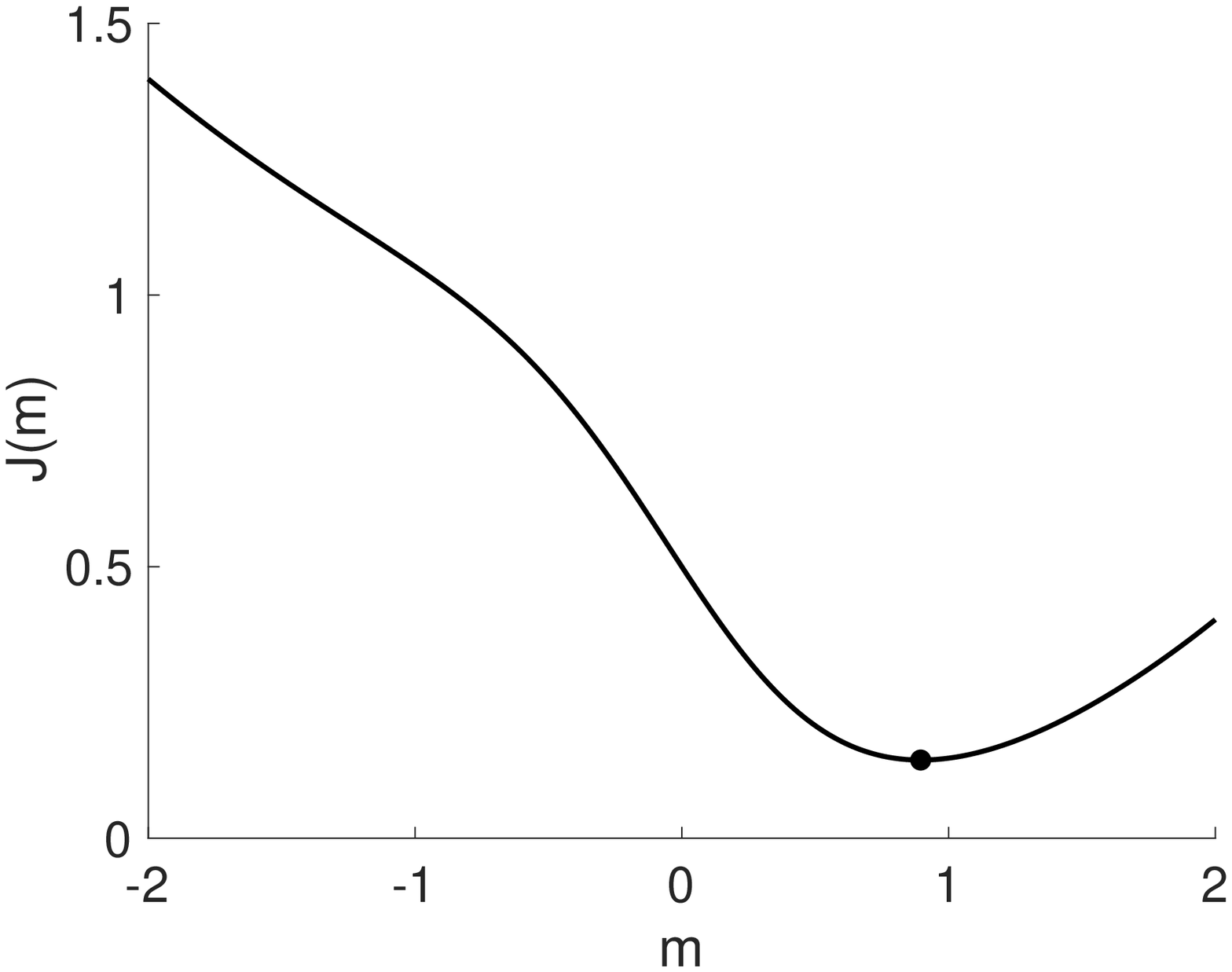}
\includegraphics[width=0.32\textwidth]{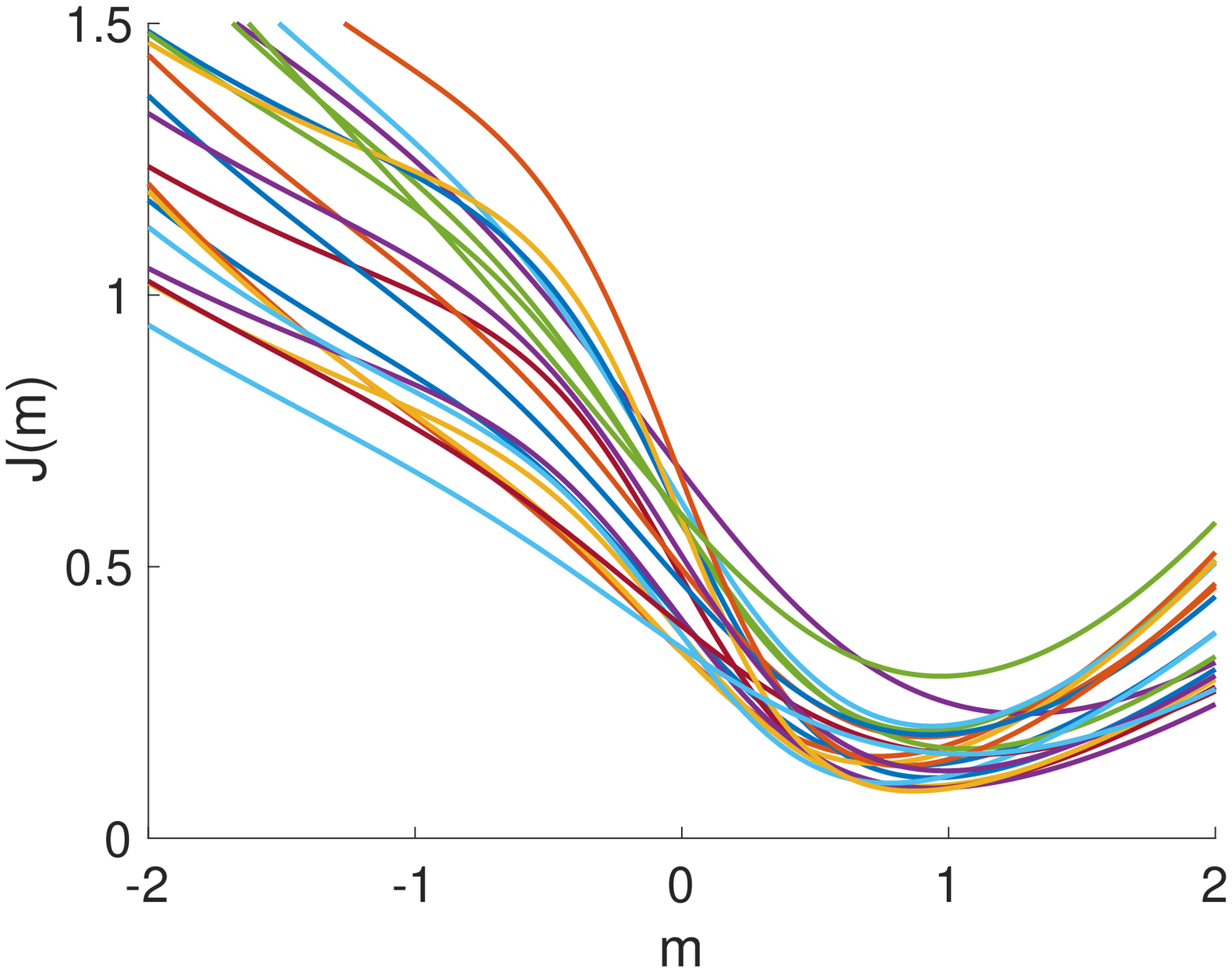}
\caption{The nominal model $J(m, \bar{\vec\theta})$~(left) and
several realizations of the model~(right) for the one-dimensional example. 
In the left panel, the black dot indicates 
the location of the minimizer.}
\label{fig:rlz}
\end{figure}

Consider the function
\[
    J(m, \vec\theta) = \frac{\theta_1}{1 + e^{\theta_2 m}} 
    + \theta_3 m^2, \quad m \in \R, \vec\theta \in \Theta.
\]
Here 
$\Theta \subset\R^3$ corresponds to taking $40\%$ of $\vec\theta$ around their nominal values $\bar{\vec\theta} = [\begin{matrix} 1 & 3 & 0.1\end{matrix}]^\top$. In Figure~\ref{fig:rlz}~(left), we show 
$J$ with $\vec\theta = \bar{\vec\theta}$, and 
in Figure~\ref{fig:rlz}~(right), we display several realizations of the function $J(m, 
\vec\theta)$, corresponding to random draws from the (uniform) distribution of $\vec\theta$. 
This demonstrates significant variations in the location of the minimizer.

\begin{figure}[ht]\centering
\includegraphics[width=0.32\textwidth]{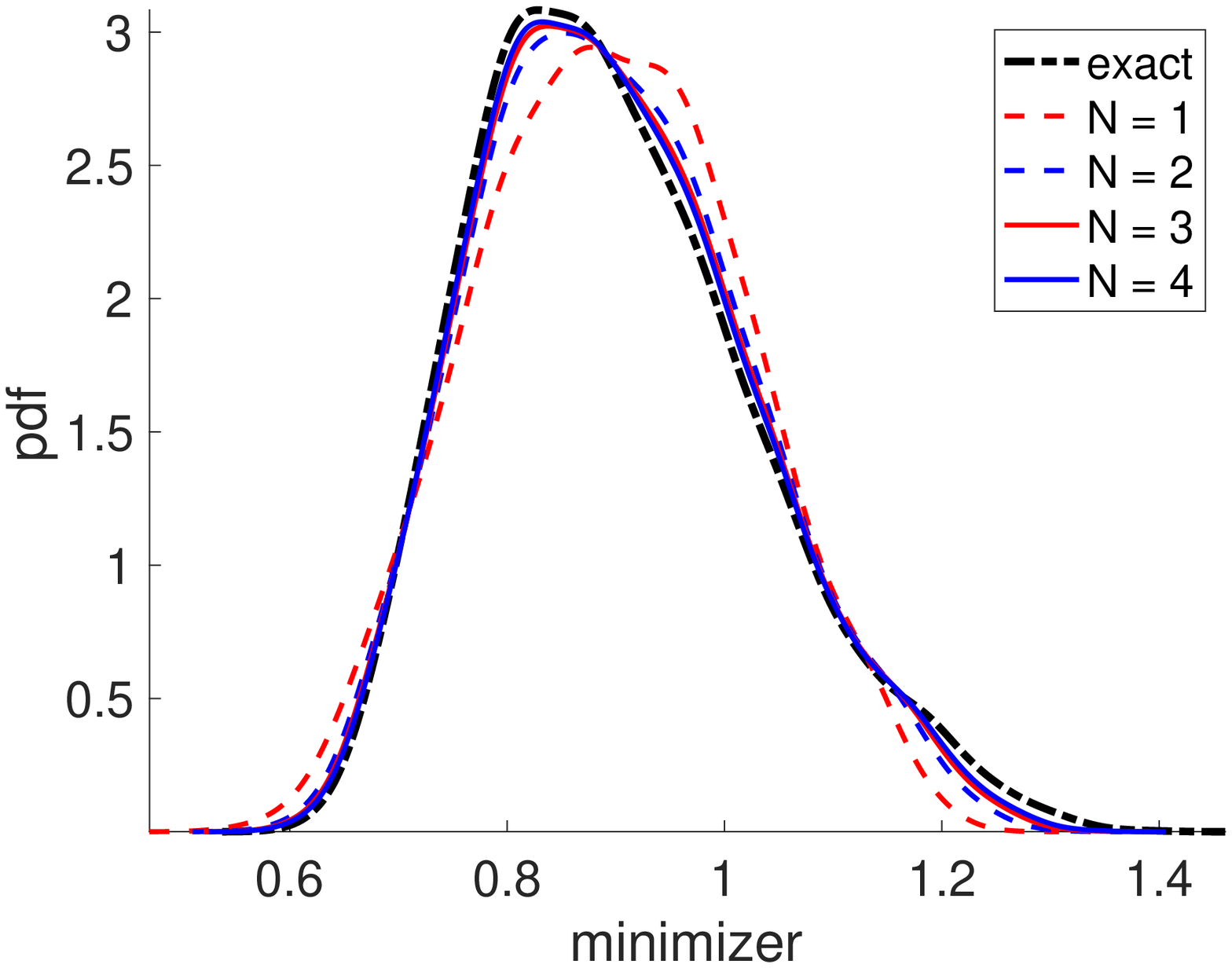}
\includegraphics[width=0.32\textwidth]{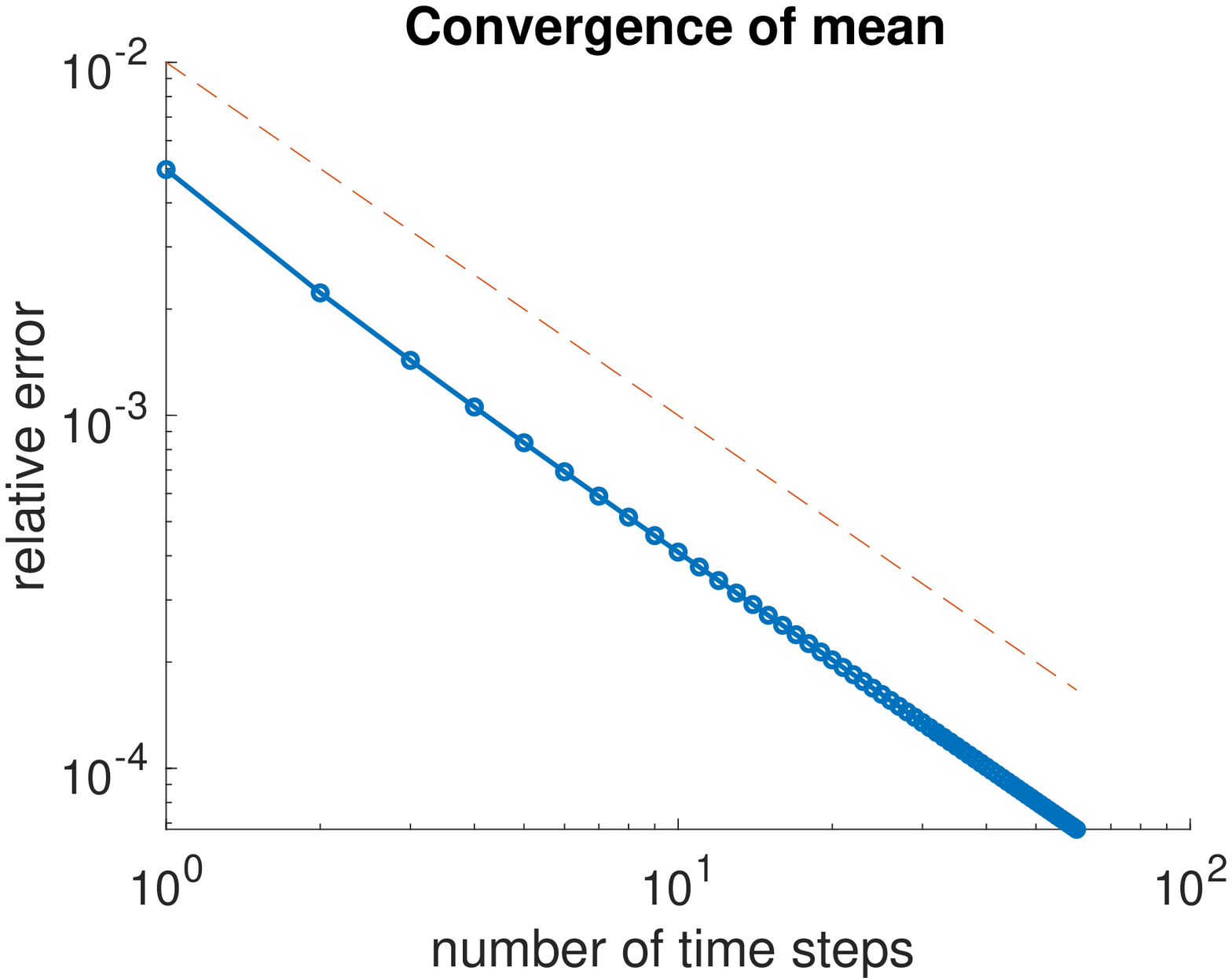}
\includegraphics[width=0.32\textwidth]{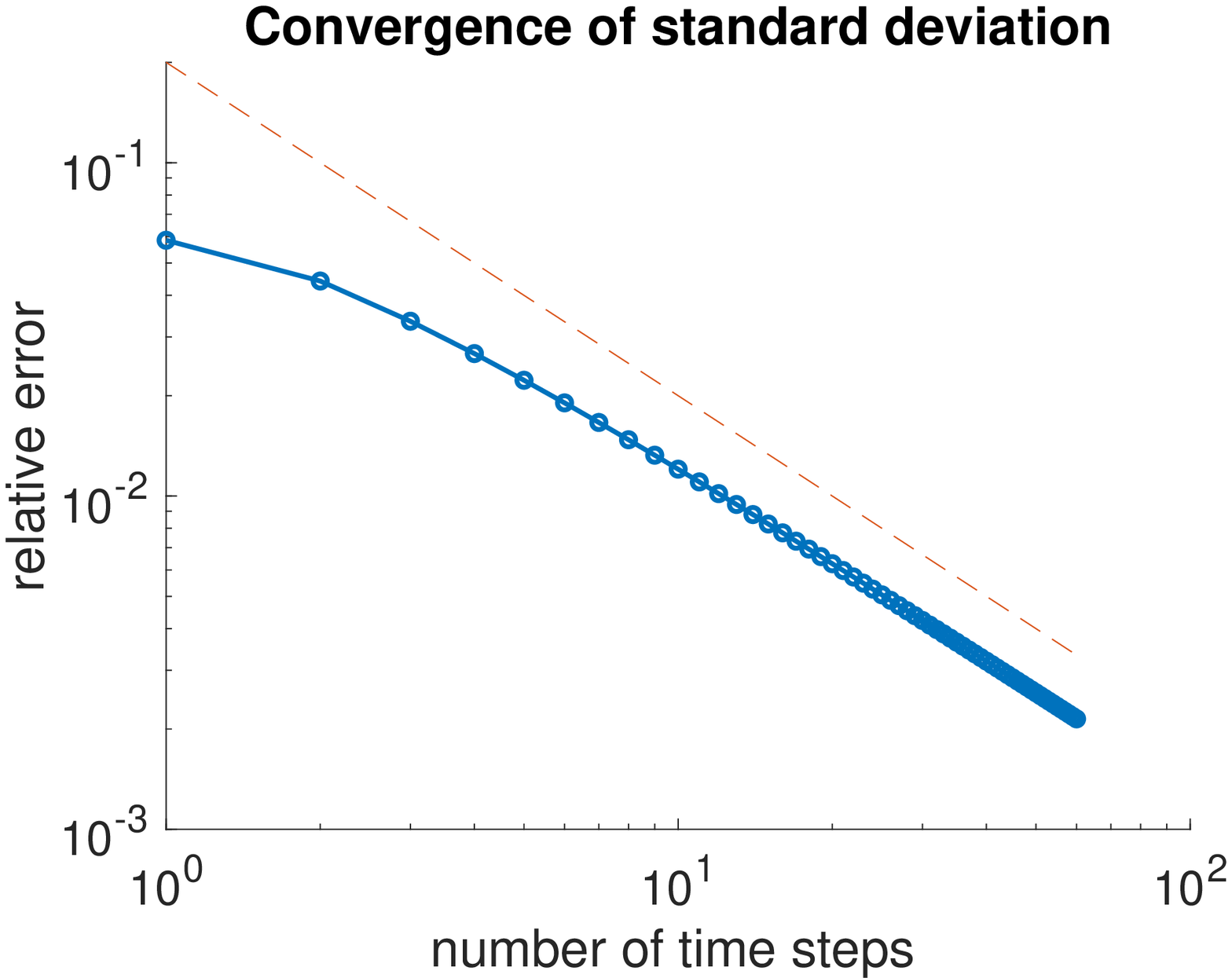}
\caption{Convergence of the pdfs~(left), 
the mean~(middle), and standard deviation~(right), 
as $N$ increases, for the one-dimensional example. In the middle and
right panels, the dashed lines indicate $\mathcal{O}(h)$, where $h = 1/N$.}
\label{fig:example_one}
\end{figure}

In Figure~\ref{fig:example_one}~(left), we show the probability density
function (pdf) of $\ms_N$, for a few choices of $N$ in~\eqref{equ:euler}.  We
also track the convergence of the mean and standard deviation in
Figure~\ref{fig:example_one}~(middle/right). The optimization
problem was solved for $5000$ realizations of $\vec\theta$ to generate a reference distribution for $\ms(\ta)$. 
From Figure~\ref{fig:example_one} we see that a small $N$ is
sufficient for approximating the pdf of $\ms$. In fact, $N = 1$
provides a reasonable approximation, and as $N$ grows, the pdf of $\ms_N$
approaches that of $\ms$ rapidly. Also, the mean and standard
deviation exhibit a first order convergence consistent with the convergence
rate of  forward Euler.

\subsection{A differential equation constrained example}

\begin{figure}\centering
\includegraphics[width=0.19\textwidth]{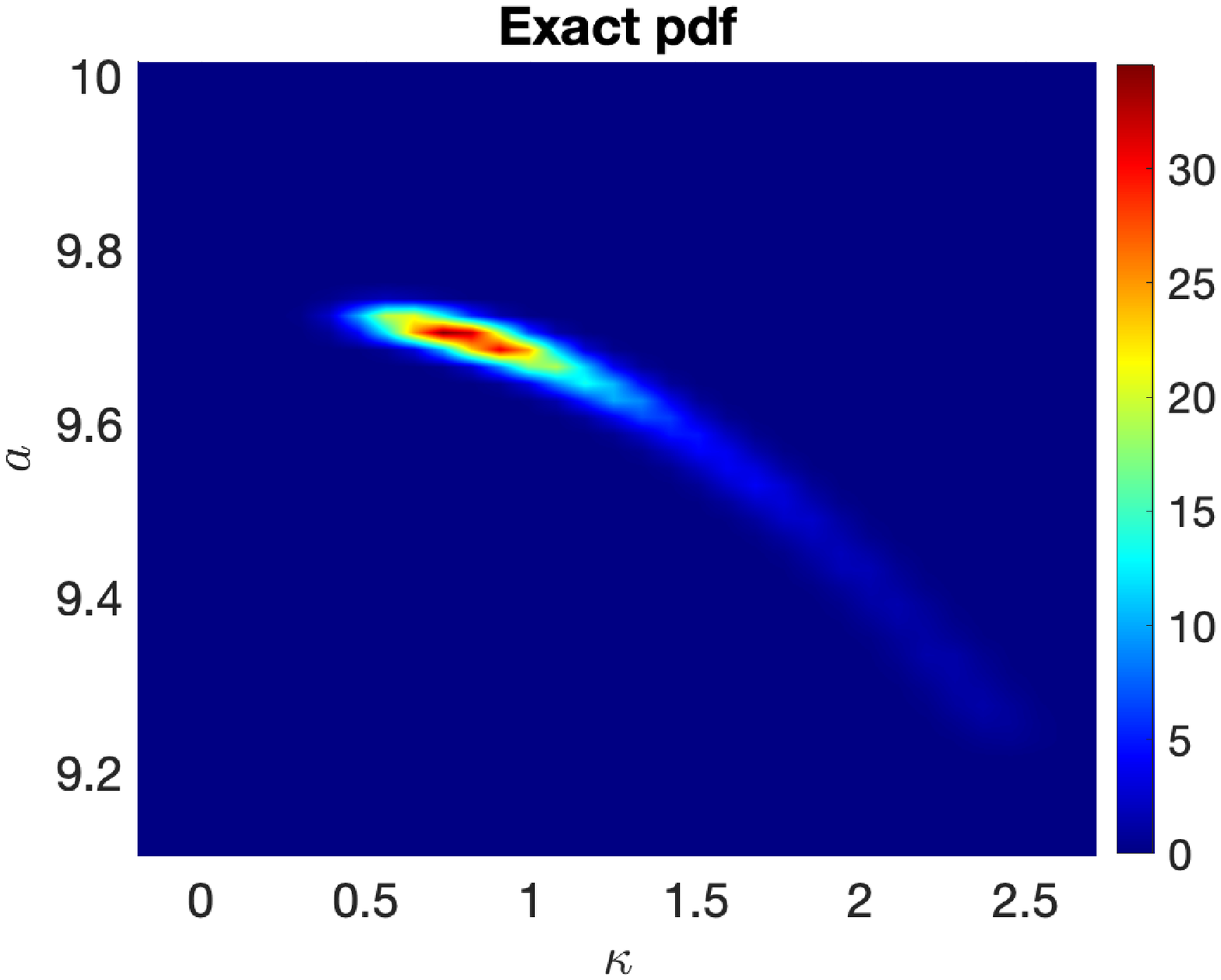}
\includegraphics[width=0.19\textwidth]{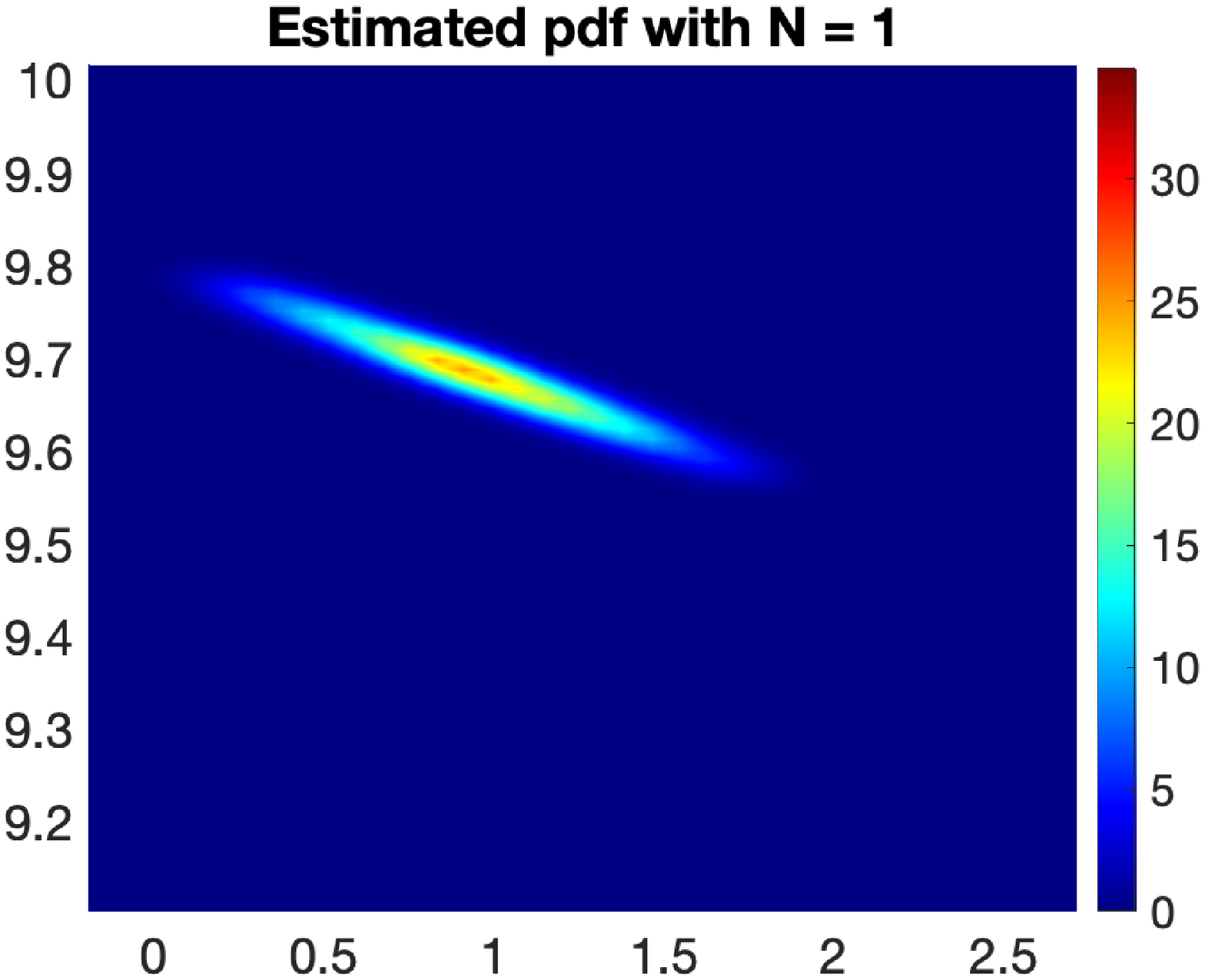}
\includegraphics[width=0.19\textwidth]{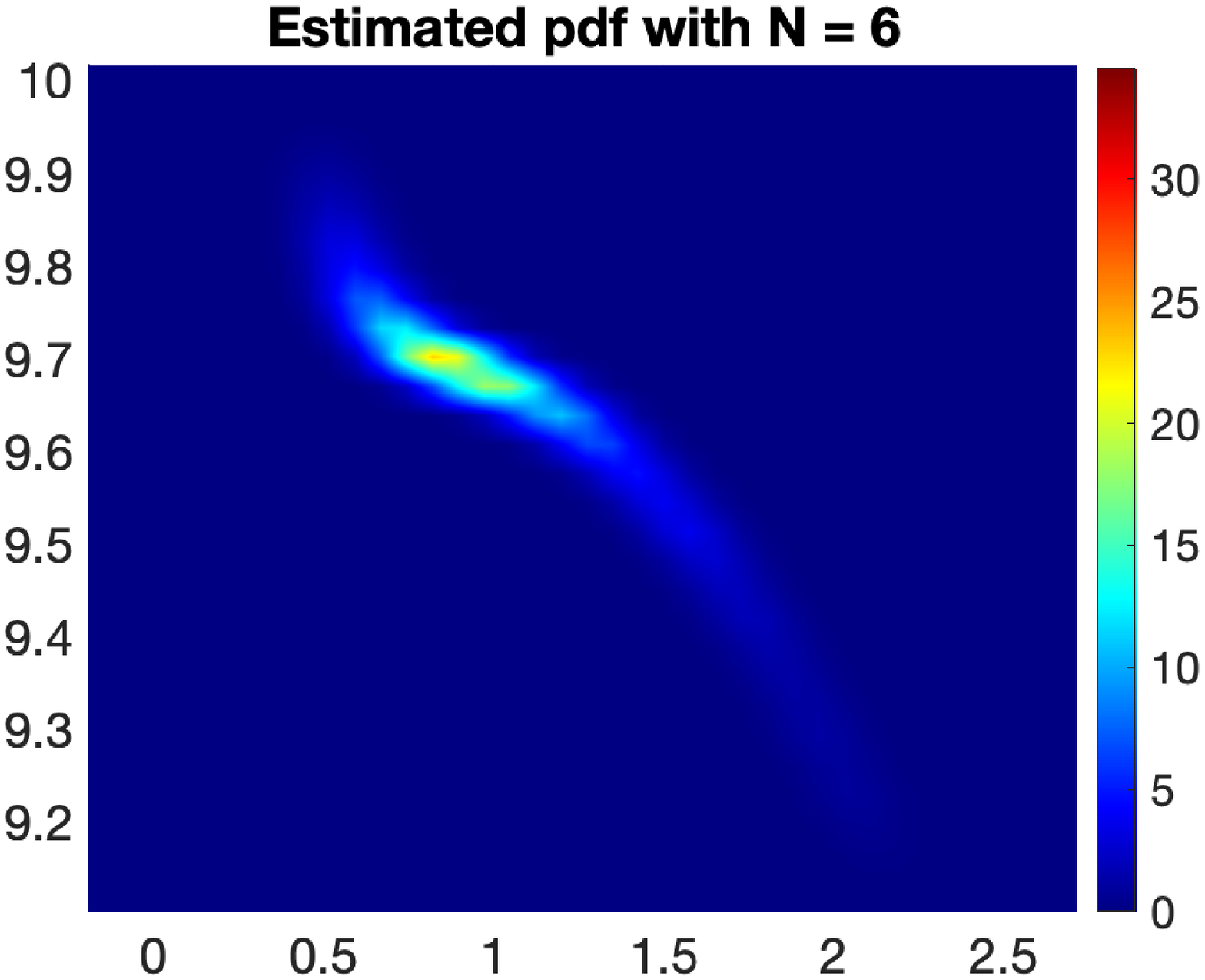}
\includegraphics[width=0.19\textwidth]{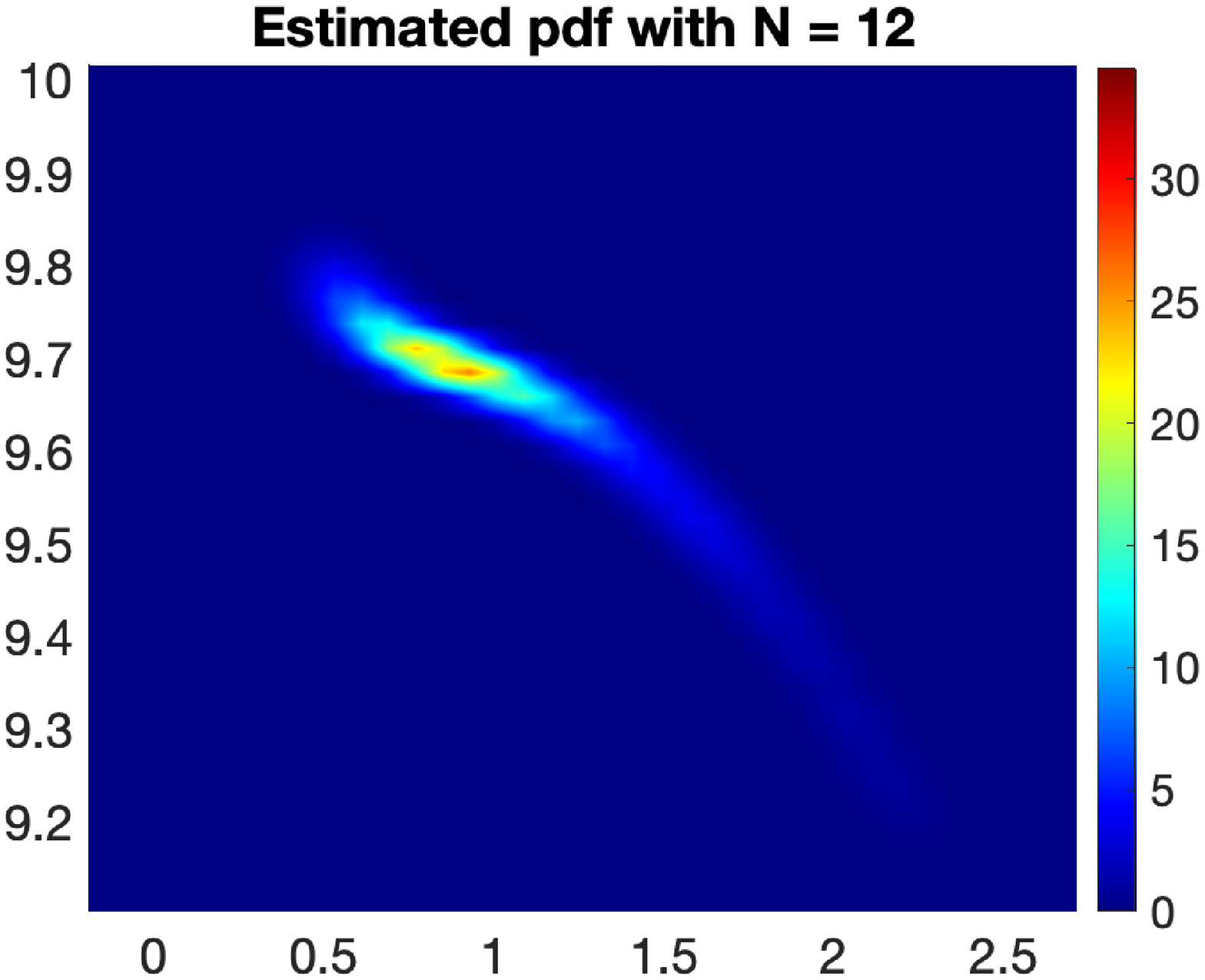}
\includegraphics[width=0.19\textwidth]{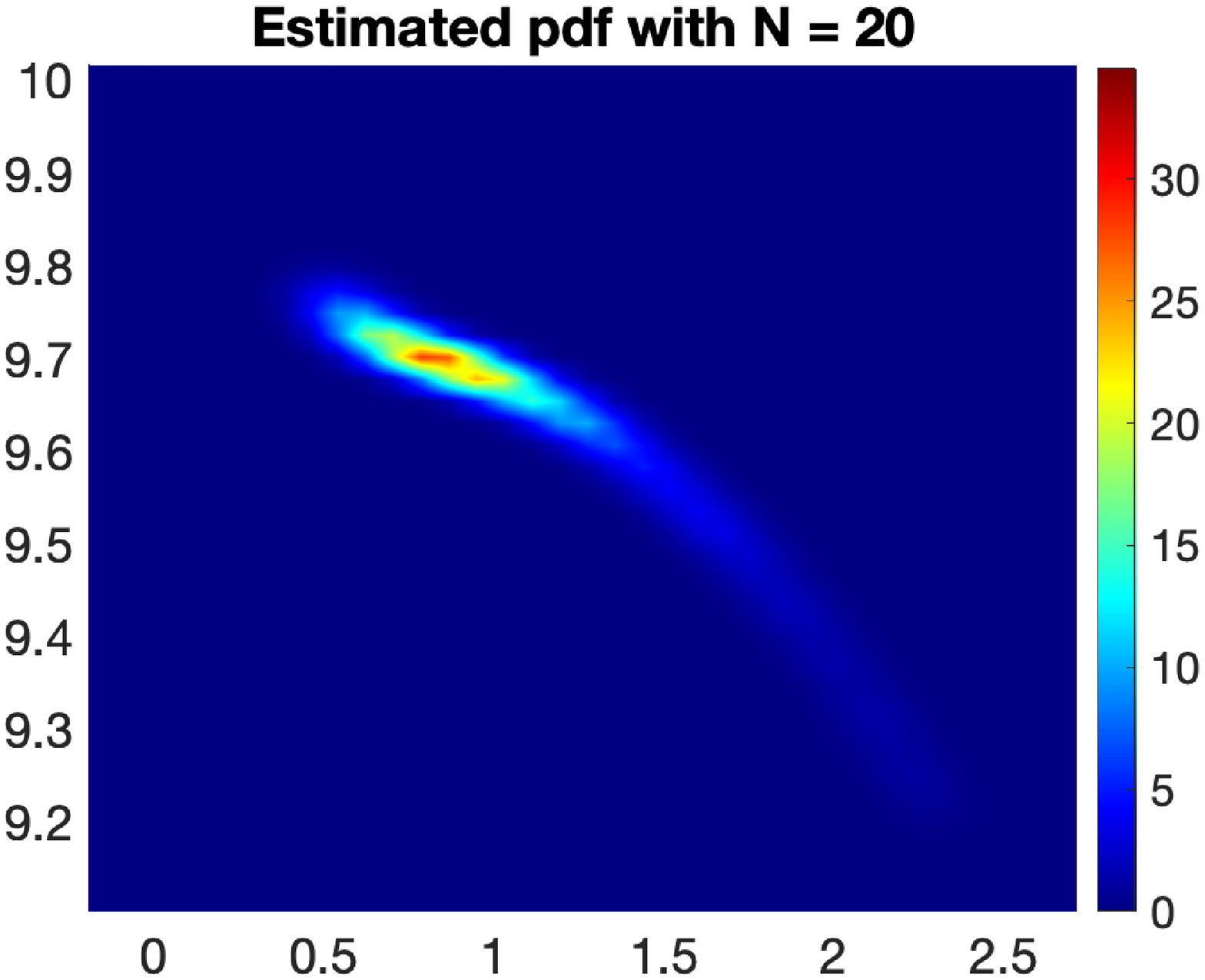}
\caption{Convergence of the joint pdfs for increasing $N$ on the differential equation constrained example.} 
\label{fig:PDF_conv_diff_eq_example}
\end{figure}

\begin{figure}\centering
\includegraphics[width=0.24\textwidth]{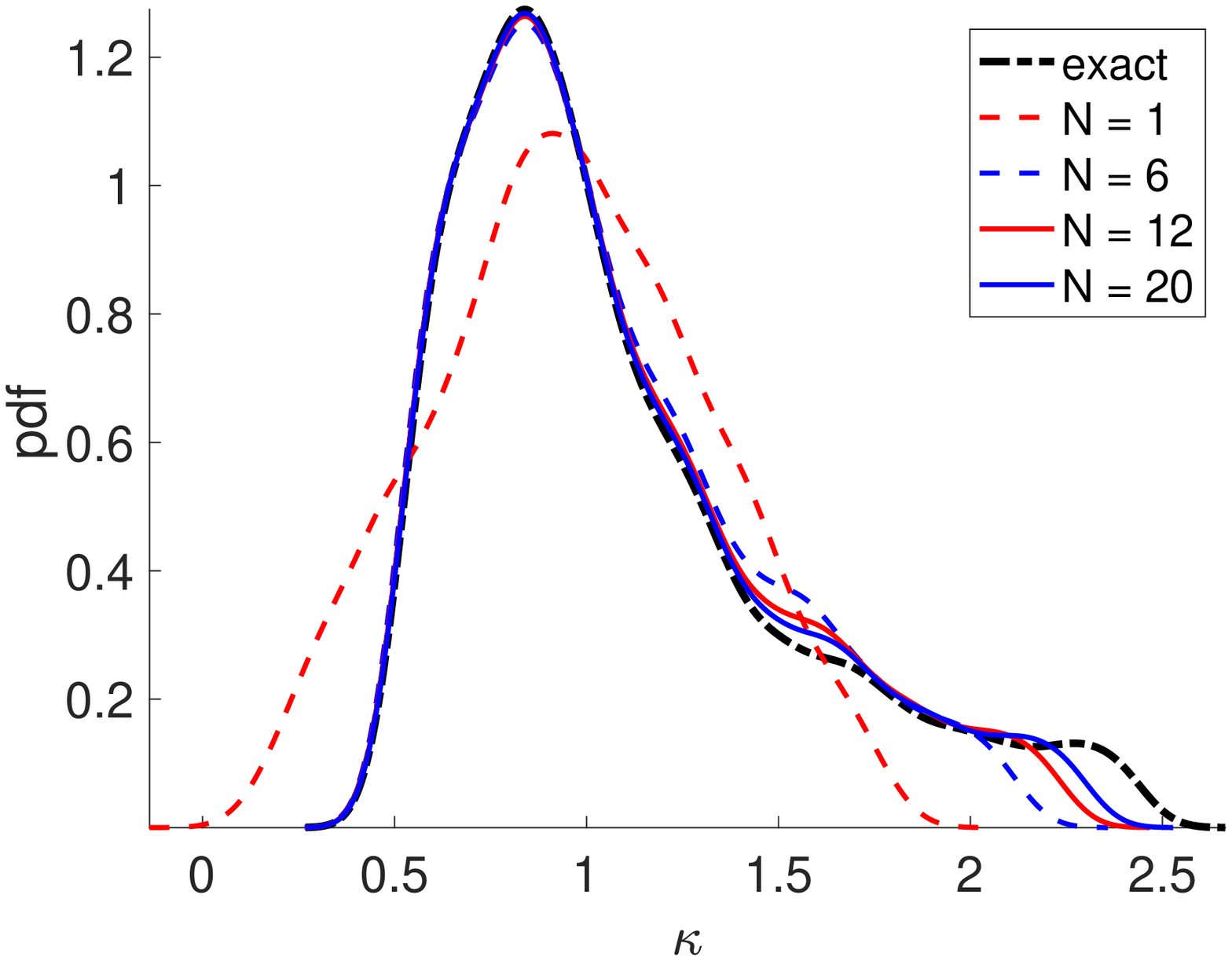}
\includegraphics[width=0.24\textwidth]{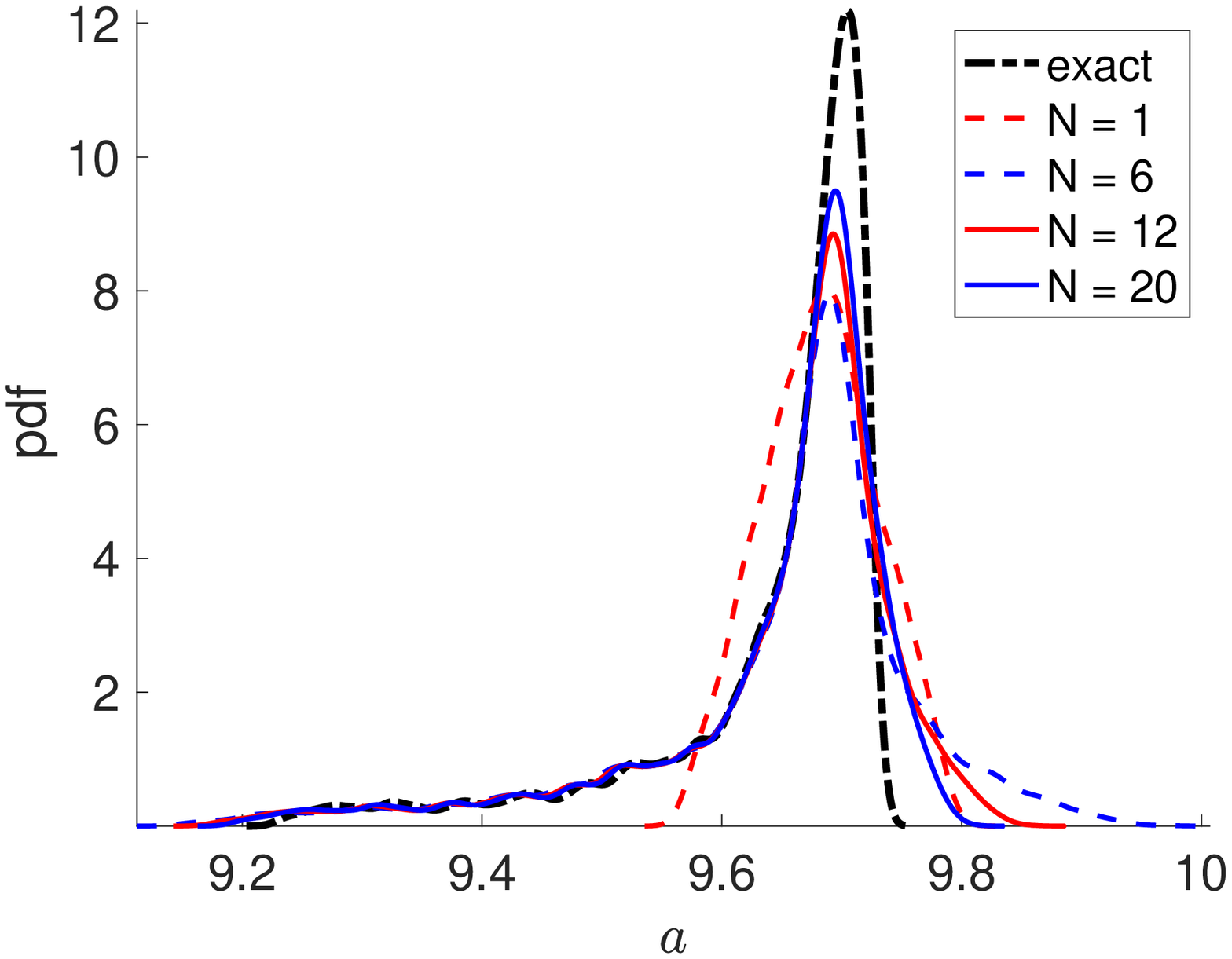}
\caption{Convergence of the marginal pdfs for increasing $N$ on the differential equation constrained example.} 
\label{fig:PDF_conv_diff_eq_example_marginals}
\end{figure}

We revisit our illustrative optimization problem~\eqref{equ:optim_CD} and
consider estimating the pdf of the optimal solution. We draw 5000 parameter
samples from a uniform distribution modeling $\pm 20 \%$ uncertainty around the
nominal parameter vector $\bar{\vec\theta} =  [\begin{matrix} 10.0 & 0.05 & 1.0
\end{matrix}]^\top$. Figure~\ref{fig:PDF_conv_diff_eq_example} displays the
joint pdf computed by solving the optimization problem for each parameter
sample (left) and compares it with the estimated joint pdf coming from our
proposed approach using $N=1,6,12$, and $20$ time steps. In each of these cases,
we solve the optimization problem once for $\vec\theta=\bar{\vec\theta}$ and
then solve the IVP~\eqref{equ:IVP} for each parameter sample to estimate the
minimizer. Using the same samples,
Figure~\ref{fig:PDF_conv_diff_eq_example_marginals} shows the convergence of
marginal pdfs of the minimizer. We observe that some information about the
correlation structure in the joint pdf and the marginal pdf for $\kappa$ can be
inferred with a small $N$. The marginal pdf of $a$ exhibits complex feature that are 
not easily resolved with a small $N$; nonetheless, for modest values of $N$ we are
able to capture many of its features. 

\section{Conclusion}\label{sec:conc}
The time-stepping approach introduced in this article is a new perspective on a
classical problem of studying parametric uncertainty in optimization problems.
Tools from post-optimality sensitivity analysis have traditionally been used to
perform local parameter studies~\cite{shapiro_SIAM_review}. By formulating an
ordinary differential equation driven by the post-optimality sensitivity, our
approach offers a mathematically rigorous approach to transition from local to
global parameter studies. We conjecture that with suitable time
discretizations, the computational cost of our approach can be less than the
cost of resolving the optimization problem for each new parameter sample. Many
questions remain.  These include the trade-offs between higher order time
integration schemes, finer temporal discretizations, and stability of the time
stepping. Further, there may be opportunities to leverage information from
previous time steps for preconditioning of future solves or reusing time steps
to explore multiple parameter samples. Exploiting such structure may enable
further gains in the computational performance of our approach relative to the
base line of repeatedly resolving optimization problems. 

Another area of inquiry is to extract global sensitivity information alongside
the distribution of the optimal solution. Each time step computes the action of
the post-optimality sensitivity operator for a different sample. Understanding
how to aggregate this sensitivity information will provide valuable insights
that are not available from repeated optimization solves. Such global aggregation 
of sensitivity information is common in the derivative-based global sensitivity
analysis literature~\cite{dgsm1,dgsm2} and uncertainty quantification more 
broadly~\cite{Smith13}.

\section*{Acknowledgments}

This article has been authored by an employee of National Technology \& Engineering Solutions of Sandia, LLC under Contract No. DE-NA0003525 with the U.S. Department of Energy (DOE). The employee owns all right, title and interest in and to the article and is solely responsible for its contents. The United States Government retains and the publisher, by accepting the article for publication, acknowledges that the United States Government retains a non-exclusive, paid-up, irrevocable, world-wide license to publish or reproduce the published form of this article or allow others to do so, for United States Government purposes. The DOE will provide public access to these results of federally sponsored research in accordance with the DOE Public Access Plan https://www.energy.gov/downloads/doe-public-access-plan. 

Sandia National Laboratories is a multimission laboratory managed and operated by National Technology and Engineering Solutions of Sandia, LLC., a wholly owned subsidiary of Honeywell International, Inc., for the U.S. Department of Energy’s National Nuclear Security Administration under contract DE-NA0003525. This paper describes objective technical results and analysis. Any subjective views or opinions that might be expressed in the paper do not necessarily represent the views of the U.S. Department of Energy or the United States Government. SAND2022-12939 O. 

This work was supported by the US Department of Energy, Office of Advanced Scientific Computing Research, Field Work Proposal 20-023231.

The work of A.~Alexanderian and M.~Stevens was supported in part by National
Science Foundation under grant DMS-1745654. The work of A.~Alexanderian was
also supported in part by the National Science Foundation under grant
DMS-2111044.

\bibliographystyle{elsarticle-num.bst} 

\bibliography{refs}

\end{document}